\documentclass[11pt]{article}
\usepackage{amsthm}
\usepackage{pxfonts}
\usepackage{yfonts}
\usepackage{dsfont}
\usepackage{graphicx}
\usepackage{relsize}
\usepackage{color}

\def\bel{\begin{equation}\label}
\def\eeq{\end{equation}}
\def\ds{\displaystyle}
\def\endproof{\hphantom{MM}
\hfill\llap{$\square$}\goodbreak}

\def\mt{\longrightarrow}
\def\v{\vskip 1em}
\def\vsk{\vskip 40em}
\def\ve{\varepsilon}
\def\R{\mathbb R}

\def\C{\mathfrak{C}}

\def\N{{\bf N}}

\def\P{{\bf P}}
\def\Q{{\bf Q}}

\def\A{{\bf A}}
\def\B{{\bf B}}

\def\L{{\bf L}}
\def\U{{\bf U}}
\def\V{{\bf V}}
\def\T{{\bf T}}

\def\I{{\bf I}}
\def\II{{\bf II}}

\def\alpha{\alphaup}
\def\beta{\betaup}
\def\gamma{\gammaup}
\def\delta{\deltaup}
\def\xi{{\xiup}}
\def\eta{{\etaup}}
\def\tau{{\tauup}}
\def\rho{{\rhoup}}
\def\phi{{\phiup}}
\def\psi{{\psiup}}
\def\lambda{{\lambdaup}}
\def\omega{\omegaup}
\def\varphi{{\varphiup}}
\def\gamma{{\gammaup}}

\def\t{{\bf t}}

\def\vol{{\bf vol}}

\newtheorem{remark}{Remark}[section]

\begin{document}
\[\hbox{\LARGE{\bf On the end-point of Stein-Weiss inequality}}\]

\[\hbox{Chuhan Sun~~~~~and~~~~~~Zipeng Wang}\]
\begin{abstract}
This paper has two purposes. First, we show that the classical 
Stein-Weiss inequality is true for $p=1$. Second, by considering a family of  strong fractional integral operators whose kernels have singularity  on every coordinate subspace, we extend this end-point result to the multi-parameter settings.
\end{abstract}

\section{Introduction}
 \setcounter{equation}{0}
 Let $0<\alpha<\N$ and $\gamma,\delta<\N$. We define
\bel{i_alpha}
\I_{\alpha\gamma\delta} f(x)~=~\int_{\R^\N} f(y)\left({1\over |x|}\right)^\gamma\left({1\over |x-y|}\right)^{\N-\alpha} \left({1\over |y|}\right)^\delta dy,\qquad \hbox{\small{$x\neq0$}}.
\eeq
In 1928, Hardy and Littlewood \cite{Hardy-Littlewood} first investigated the $\L^p-\L^q$-regularity of $\I_{\alpha\gamma\delta}$ at $\N=1$. Thirty years later,  Stein and Weiss \cite{Stein-Weiss} extended this result to every higher dimensional space.  Today, it bears the name of Stein-Weiss inequality.  

$\diamond$ {\small Throughout, $\C>0$ is regarded as a generic constant depending on its subindices}.
\v
{\bf Theorem A:~ Stein and Weiss, 1958}  

{\it Let $\I_{\alpha\gamma\delta}$ defined in (\ref{i_alpha}) for $0<\alpha<\N$ and $\gamma,\delta<\N$. We have
\bel{norm ineq}
\left\|\I_{\alpha\gamma\delta} f\right\|_{\L^q\left(\R^\N\right)}~\leq~\C_{\alpha~\gamma~\delta~p~q}~\left\| f\right\|_{\L^p\left(\R^\N\right)},\qquad 1<p\leq q<\infty
\eeq
 if and only if
\bel{constraints}
\gamma~<~{\N\over q},\qquad \delta~<~\N\left({p-1\over p}\right),\qquad \gamma+\delta~\ge~0,\qquad 
{\alpha\over \N}~=~{1\over p}-{1\over q}+{\gamma+\delta\over \N}.
\eeq}
\begin{remark}
In the original  paper of Stein and Weiss \cite{Stein-Weiss}, (\ref{constraints}) is given as a sufficient condition. These constraints of $\alpha,\gamma,\delta,p,q$ in (\ref{constraints}) are in fact necessary. See subsection {\bf 3.1}.
\end{remark}
When $\gamma=\delta=0$, {\bf Theorem A}  was proved in $\R^\N$ by Sobolev \cite{Sobolev}. This is also known   as Hardy-Littlewood-Sobolev inequality. 

The theory of  fractional integration in  weighted norms has been substantially developed during the second half of  20th century. See   Coifman and Fefferman \cite{Coifman-C.Fefferman}, Fefferman and Muckenhoupt \cite{C.Fefferman-Muckenhoupt},   Muckenhoupt and Wheeden \cite{Muckenhoupt-Wheeden},   P\'{e}rez \cite{Perez} and Sawyer and Wheeden \cite{Sawyer-Wheeden}. 

Our first main result  is an improvement of {\bf Theorem A} to include the end-point $p=1$. 

{\bf Theorem One}  ~~{\it Let  $\I_{\alpha\gamma\delta}$ defined in (\ref{i_alpha}) for $0<\alpha<\N$ and $\gamma, \delta<\N$. We have
\bel{norm ineq L^1}
\left\|\I_{\alpha\gamma\delta} f\right\|_{\L^q\left(\R^\N\right)}~\leq~\C_{\alpha~\gamma~\delta~p~q}~\left\| f\right\|_{\L^p\left(\R^\N\right)},\qquad 1\leq p\leq q<\infty
\eeq
 if and only if
\bel{constraints+formula}
\gamma~<~{\N\over q},\qquad \delta~<~\N\left({p-1\over p}\right),\qquad \gamma+\delta~\ge~0,\qquad {\alpha\over \N}~=~{1\over p}-{1\over q}+{\gamma+\delta\over\N}.
\eeq}
\begin{remark}
For $p=1$ and $\gamma+\delta>0$,   {\bf Theorem One} is recently proved by N\'{a}poli  and Picon \cite{Napoli-Picon}.
  The 'good kernel' estimate used thereby is nice but cannot handle the case when $\gamma+\delta=0$. 
\end{remark}
To obtain the {\it if and only if} statement in (\ref{norm ineq L^1})-(\ref{constraints+formula}) when $p=1$,  we need to go through  an argument of interpolation with changing measures.
\v

Next, consider $\R^\N=\R^{\N_1}\times\R^{\N_2}\times\cdots\times\R^{\N_n}$.
Let $\gamma, \delta<\N$ and
\bel{alpha_i s}
\alpha~=~\alpha_1+\alpha_2+\cdots+\alpha_n,\qquad 0<\alpha_i<\N_i,\qquad  i=1,2,\ldots,n.
\eeq 
We define
\bel{I_alpha f}
\II_{\alpha\gamma\delta} f(x)~=~\int_{\R^\N} f(y)\left({1\over |x|}\right)^\gamma \prod_{i=1}^n \left({1\over |x_i-y_i|}\right)^{\N_i-\alpha_i} \left({1\over |y|}\right)^\delta dy,\qquad \hbox{\small{$x\neq0$}}.
\eeq
The study of  certain operators  that  commute with a multi-parameter family of dilations  dates back to the time of Jessen, Marcinkiewicz and Zygmund. Over the several past decades, a number of remarkable  results have been accomplished, for example 
by Fefferman \cite{R.Fefferman}, C\'{o}rdoba and Fefferman \cite{Cordoba-Fefferman},  Fefferman and Stein \cite{R.Fefferman-Stein},  M\"{u}ller, Ricci and Stein \cite{M.R.S}, Journ\'{e} \cite{Journe'} and Pipher \cite{Pipher}.   

For  $1<p\leq q<\infty$,  a characterization is established  between  the norm inequality 
\[\left\|\II_{\alpha\gamma\delta} f\right\|_{\L^q\left(\R^\N\right)}~\leq~\C_{p~q~\alpha~\gamma~\delta}\left\| f\right\|_{\L^p\left(\R^\N\right)},\qquad 1<p\leq q<\infty\]
and the necessary constraints  of $\gamma, \delta,p,q$ and $\alpha_i, i=1,2,\ldots,n$. 
See another recent work  \cite{Wang}.
This extends {\bf Theorem A} to the multi-parameter settings.
Our second result is to give such an extension when $p=1$.    
\v
{\bf Theorem Two}  ~~{\it Let  $\II_{\alpha\gamma\delta}$ defined in (\ref{I_alpha f}) for $0<\alpha_i<\N_i, i=1,2,\ldots,n$ and $\gamma, \delta<\N$. We have
\bel{Norm Ineq L^1}
\left\|\II_{\alpha\gamma\delta} f\right\|_{\L^q\left(\R^\N\right)}~\leq~\C_{\alpha~\gamma~\delta~q}~\left\| f\right\|_{\L^1\left(\R^\N\right)},\qquad 1\leq q<\infty
\eeq
 if and only if
\bel{Formula}
\gamma~<~{\N\over q},\qquad \delta~<~0,\qquad \gamma+\delta~\ge~0,\qquad  {\alpha\over\N}~=~1-{1\over q}+{\gamma+\delta\over\N}
\eeq
and
\bel{Constraints}
\alpha_i-\N_i~<~\delta,\qquad i~=~1,2,\ldots,n.
\eeq}

In order to prove {\bf Theorem Two}, we develop a new framework where the product space $\R^{\N_1}\times\R^{\N_2}\times\cdots\times\R^{\N_n}$ is decomposed into an infinitely many dyadic cones. Every  partial operator  defined on one of these  cones is essentially an  one-parameter fractional integral operator, satisfying the regarding $\L^1\mt\L^q$-norm inequality. Furthermore, the operator's norm decays exponentially as the eccentricity of the cone getting large.

This paper is organized as follows. The next section is devoted to two fundamental lemmas initially given by Stein and Weiss \cite{Stein-Weiss} for $p>1$. These results are improved now to become applicable for $p\ge1$. We prove
{\bf Theorem One} and {\bf Theorem Two} in section 3 and  section 4 respectively.

\section{Two fundamental lemmas}
\setcounter{equation}{0}
{\bf Lemma One}~~{\it Let $\Omega(u,v)\geq0$  defined in the quadrant $\left\{(u,v):u\geq0,~ v\geq0\right\}$ which is homogeneous of degree $-\N$ and
\bel{J}
\A~\doteq~\int_{0}^{\infty}\Omega(1,t)t^{\N\big({p-1\over p}\big)-1}dt~<~\infty.
\eeq
Consider 
\bel{U f}
\U f(x)~=~\int_{\R^\N} \Omega\left(|x|,|y|\right)f(y)dy.
\eeq
We have
\bel{Ua Ineq}
\left\|\U f\right\|_{\L^p(\R^\N)}~\leq~\C_\A~\|f\|_{\L^p(\R^\N)},\qquad 1\leq p<\infty.
\eeq
}\\
{\bf Proof}~~Let $R=|x|$ and $r=|y|$. For $\N\ge2$, write $x=R\xi$ and $y=r\eta$ of which $\xi,\eta$ are  unit vectors.
We have 
\bel{UR}
\U f(x)~=~\int_{\mathds{S}^{\N-1}}\int_0^\infty \Omega(R,r)f(r\eta)r^{\N-1}drd\sigma(\eta)
\eeq
where $\sigma$ denotes the surface measure on $\mathds{S}^{\N-1}$. 

Consider
\bel{U eta R Ineq}
\begin{array}{lr}\ds
\left\{\int_0^{\infty}\left|\int_0^\infty \Omega(R,r)f(r\eta)r^{\N-1}dr\right|^pR^{\N-1}dR\right\}^{1\over p}
\\\\ \ds
~=~\left\{\int_0^{\infty}\left|\int_0^\infty \Omega(1,t)f(tR\eta)t^{\N-1}dt\right|^p R^{\N-1}dR\right\}^{1\over p}
\qquad
 \hbox{\small{$r=tR$ and  $\Omega$ is homogeneous of degree $-\N$}}
\\\\ \ds
~\leq~\int_0^\infty \Omega(1,t)t^{\N-1}\left\{\int_0^\infty \left|f(tR\eta)\right|^pR^{\N-1}dR\right\}^{1\over p}dt
\qquad
 \hbox{\small{ by Minkowski integral inequality}}
\\\\ \ds
~=~\int_0^\infty \Omega(1,t)t^{\N\big[1-{1\over p}\big]-1}\left\{\int_0^\infty \left|f(r\eta)\right|^pr^{\N-1}dr\right\}^{1\over p}dt
\\\\ \ds
~=~\A~\left\{\int_0^{\infty}\left|f(r\eta)\right|^pr^{\N-1}dr\right\}^{1\over p}.
\end{array}
\eeq
We find
\bel{Ua LSH}
\begin{array}{lr}\ds
\left\|\U f\right\|_{\L^p(\R^\N)}~=~\left\{\int_{\mathds{S}^{\N-1}} \int_0^\infty \left|\int_{\mathds{S}^{\N-1}}\int_0^\infty \Omega(R,r)f(r\eta)r^{\N-1}drd\sigma(\eta)\right|^p R^{\N-1}dRd\sigma(\xi)\right\}^{1\over p}
\\\\ \ds~~~~~~~~~~~~~~~~~~~
~=~\omega_{\N-1}^{1\over p} \left\{\int_0^\infty \left|\int_{\mathds{S}^{\N-1}}\int_0^\infty \Omega(R,r)f(r\eta)r^{\N-1}drd\sigma(\eta)\right|^p R^{\N-1}dR\right\}^{1\over p}
\\\\ \ds
~~~~~~~~~~~~~~~~~~~\leq~\omega_{\N-1}^{1\over p} \int_{\mathds{S}^{\N-1}}\left\{\int_0^\infty\left|\int_0^\infty \Omega(R,r)f(r\eta)r^{\N-1}dr\right|^p R^{\N-1}dR \right\}^{1\over p} d\sigma(\eta) 
\\
~~~~~~~~~~~~~~~~~~~~~~~~~~~~~~~~~~~~~~~~~~~~~~~~~~~~~~~~~~~~~~~~~~~~~~\hbox{\small{ by Minkowski integral inequality}}
\\ \ds
~~~~~~~~~~~~~~~~~~
~\leq~\omega_{\N-1}^{1\over p} \A \int_{\mathds{S}^{\N-1}}\left\{\int_0^{\infty}\left|f(r\eta)\right|^pr^{\N-1}dr\right\}^{1\over p} d\sigma(\eta)\qquad\hbox{\small{by (\ref{U eta R Ineq})}}
\\\\ \ds
~~~~~~~~~~~~~~~~~~
~\leq~\A\omega_{\N-1}^{1\over p}\left\{  \int_{\mathds{S}^{\N-1}} \int_0^{\infty} \left|f(r\eta)\right|^pr^{\N-1}dr d\sigma(\eta)\right\}^{1\over p} \left\{ \int_{\mathds{S}^{\N-1}} d\sigma(\eta)\right\}^{p-1\over p}
\\ \ds~~~~~~~~~~~~~~~~~~~~~~~~~~~~~~~~~~~~~~~~~~~~~~~~~~~~~~~~~~~~~~~~~~~~~~~~~~~~~~~~~~~
\hbox{\small{by H\"{o}lder inequality}}
\\ \ds~~~~~~~~~~~~~~~~~
~=~ \A \omega_{\N-1}\left\| f\right\|_{\L^p(\R^\N)}
\end{array}
\eeq
where $\omega_{\N-1}=2\pi^{\N\over2}\Gamma^{-1}\left({\N\over 2}\right)$ is the  area of $\mathds{S}^{\N-1}$.

When $\N=1$, simply take $d\sigma$ to be the point measure on $1$ and $-1$. The same estimates hold in
 (\ref{U eta R Ineq})-(\ref{Ua LSH}).
\endproof
\v

{\bf Lemma Two}~~{\it 
Let $\N\ge2$. Define
$\Delta(t,\xi, \eta)= \left|1-2t\xi\cdot\eta+t^2\right|^{1\over 2}$ for $t>0$ and $\xi,\eta\in\mathds{S}^{\N-1}$. We have
\bel{integral of Delta}
\begin{array}{lr}\ds
\int_{\mathds{S}^{\N-1}}{1\over \Delta^{\N-\alpha}(t,\xi,\eta)}d\sigma(\xi)~=~\int_{\mathds{S}^{\N-1}}{1\over \Delta^{\N-\alpha}(t,\xi,\eta)}d\sigma(\eta)
\\\\ \ds~~~~~~~~~~~~~~~~~~~~~~~~~~~~~~~~~~~~~~~~
~\leq~\C~\left|1-t\right|^{-{\N-\alpha\over \N}},\qquad t\neq1,\qquad \xi,\eta\in\mathds{S}^{\N-1}.
\end{array}
\eeq
}

{\bf Proof }~~Observe that $\Delta(t,\xi,\eta)$ is symmetric $w.r.t$ $\xi$ and $\eta$. For $0<t<1$, we have
\bel{Possion}
\P(\xi, t\eta)~=~\frac{1-|t\eta|^2}{|\xi-t\eta|^\N}~=~\frac{1-t^2}{\Delta^\N(t,\xi,\eta)},\qquad \xi,\eta\in\mathds{S}^{\N-1}
\eeq
which is the Poisson kernel on the unit sphere  $\mathds{S}^{\N-1}$. 
A direct computation shows
\bel{harmonic}
\Delta_\eta \P(\xi, t\eta)~=~0,\qquad \xi\in\mathds{S}^{\N-1}
\eeq
where $\Delta_\eta$ is the Laplacian operator $w.r.t~\eta$.

By using the mean value property of harmonic functions, we find
\bel{P=1}
1~=~\P(\xi, 0)~=~\frac{1}{\omega_{\N-1}}\int_{\mathds{S}^{\N-1}}\P(\xi,t\eta)d\sigma(\eta),\qquad 0<t<1.
\eeq
This further implies
\bel{t<1}
{1\over \omega_{\N-1}}\int_{\mathds{S}^{\N-1}}\frac{1-t^2}{\Delta^\N(t,\xi,\eta)} d\sigma(\eta)~=~1,\qquad 0< t<1.
\eeq
On the other hand, write $0<s=t^{-1}<1$ for $t>1$. From (\ref{Possion}), we have
\bel{P t>1}
\begin{array}{lr}\ds
\frac{1-t^2}{\Delta^\N(t,\xi,\eta)}~=~\frac{1-t^2}{\left|1+2t\xi\cdot\eta+t^2\right|^{\frac{\N}{2}}}
~=~{t^2(s^2-1)\over t^\N  \left|s^2+2s\xi\cdot\eta+1\right|^{\frac{\N}{2}}}
\\\\ \ds~~~~~~~~~~~~~~~~~~
~=~-t^{2-\N} \frac{1-s^2}{\Delta^\N(s,\xi,\eta)}~=~-t^{2-\N} \P(\xi,s\eta),\qquad \xi,\eta\in\mathds{S}^{\N-1}.
\end{array}
\eeq
By using (\ref{P=1}) and (\ref{P t>1}), we find
\bel{t>1}
\begin{array}{lr}\ds
{1\over \omega_{\N-1}}\int_{\mathds{S}^{\N-1}}\frac{1-t^2}{\Delta^\N(t,\xi,\eta)} d\sigma(\eta)
~=~-t^{2-\N}~\frac{1}{\omega_{\N-1}}\int_{\mathds{S}^{\N-1}}\P(\xi,s\eta)d\sigma(\eta)
\\\\ \ds~~~~~~~~~~~~~~~~~~~~~~~~~~~~~~~~~~~~~~~~~~~~~~~
~=~-t^{2-\N},\qquad t>1.
\end{array}
\eeq
By putting together (\ref{t<1}) and (\ref{t>1}), we obtain
\bel{t est}
\begin{array}{lr}\ds
{1\over \omega_{\N-1}}\int_{\mathds{S}^{\N-1}}\frac{1}{\Delta^\N(t,\xi,\eta)} d\sigma(\eta)~\leq~{1\over |1-t^2|}~<~{1\over |1-t|}.
\end{array}
\eeq
Lastly, by applying H\"older's inequality, we have
\bel{Holder}
\begin{array}{lr}\ds
\int_{\mathds{S}^{\N-1}}\frac{1}{\Delta^{\N-\alpha}(t,\xi,\eta)}d\sigma(\eta)~\leq~\left\{\int_{\mathds{S}^{\N-1}}\left[\frac{1}{\Delta^{\N-\alpha}(t,\xi,\eta)}\right]^\frac{\N}{\N-\alpha}d\sigma(\eta)\right\}^{\frac{\N-\alpha}{\N}}\left\{\int_{\mathds{S}^{\N-1}}d\sigma(\eta)\right\}^{\frac{\alpha}{\N}}
\\\\ \ds~~~~~~~~~~~~~~~~~~~~~~~~~~~~~~~~~~~~~~~~
~=~\left\{\int_{\mathds{S}^{\N-1}}\frac{1}{\Delta^\N(t,\xi,\eta)}d\sigma(\eta)\right\}^{\frac{\N-\alpha}{\N}}\left(\omega_{\N-1}\right)^{\frac{\alpha}{\N}}
\\\\ \ds~~~~~~~~~~~~~~~~~~~~~~~~~~~~~~~~~~~~~~~~
~\leq~\left(\frac{1}{|1-t|}\right)^{\frac{\N-\alpha}{\N}}\left(\omega_{\N-1}\right)^{\N-\alpha\over\N}\left(\omega_{\N-1}\right)^{\frac{\alpha}{\N}}\qquad\hbox{\small{by (\ref{t est})}}
\\\\ \ds~~~~~~~~~~~~~~~~~~~~~~~~~~~~~~~~~~~~~~~~
~=~\omega_{\N-1}~|1-t|^{-\frac{\N-\alpha}{\N}}.
\end{array}
\eeq

\section{Proof of Theorem One}
\setcounter{equation}{0}
\subsection{Norm inequality in (\ref{norm ineq L^1}) implies constraints in (\ref{constraints+formula})}
{\bf Case 1}~~Let $p=1$. Denote $Q$ as any cube in $\R^\N$. Choose $f=\chi_Q$ which is an indicator function supported in $Q$. 
The norm inequality in (\ref{norm ineq L^1}) implies
\bel{characteristic}
\sup_{Q\subset\R^\N}~|Q|^{{\alpha\over \N}-1+{1\over q}}\left\{{1\over |Q|}\int_{Q} \left({1\over |x|}\right)^{\gamma q} dx\right\}^{1\over q}\left\{{1\over |Q|}\int_{Q} \left({1\over |x|}\right)^\delta dx\right\}~<~\infty.
\eeq

A standard exercise of changing dilations inside (\ref{characteristic}) shows that ${\alpha\over\N}=1-{1\over q}+{\gamma+\delta\over\N}$ is an necessary condition. Moreover, it is essential to have $\gamma<{\N\over q}$ for the local integrability of $|x|^{-\gamma q}$. 
We claim ${\alpha\over \N}-1+{1\over q}\ge0$.
Together with  ${\alpha\over\N}=1-{1\over q}+{\gamma+\delta\over\N}$, we must have $\gamma+\delta\ge0$. 

Suppose ${\alpha\over \N}-1+{1\over q}<0$. Let $Q$ centered on some  $x_o\neq0$. By shrinking $Q$ to $x_o$ and applying Lebesgue's Differentiation Theorem, we find 
\[ 
\left\{{1\over |Q|}\int_{Q} \left({1\over |x|}\right)^{\gamma q} dx\right\}^{1\over q}\left\{{1\over |Q|}\int_{Q} \left({1\over |x|}\right)^\delta dx\right\}~=~|x_o|^{-(\gamma+\delta)}~>~0.
\]
On the other hand, 
$|Q|^{{\alpha\over \N}-1+{1\over q}}\mt\infty$.
This contradicts to (\ref{characteristic}).

Let $\I_{\alpha\gamma\delta} f$ defined in (\ref{i_alpha}).  Assert $f\ge0, ~f\in\L^1(\R^\N)$ supported in the unit ball, denoted by $\B$. We have
\bel{I > example}
\begin{array}{lr}\ds
\I_{\alpha\gamma\delta}f(x)~\ge~\chi(|x|>10)\int_\B f(y) \left({1\over |x|}\right)^\gamma \left({1\over |x-y|}\right)^{\N-\alpha} \left({1\over |y|}\right)^\delta dy
\\\\ \ds~~~~~~~~~~~~~
~>~2^{\alpha-\N} \left({1\over |x|}\right)^{\N-\alpha+\gamma} \chi(|x|>10)~\int_\B f(y) \left({1\over |y|}\right)^\delta dy.
\end{array}
\eeq
Observe that if $\delta\ge0$, then ${\alpha\over\N}=1-{1\over q}+{\gamma+\delta\over\N}$ implies $\N-\alpha+\gamma\leq{\N\over q}$. Consequently, $\big[\I_{\alpha\gamma\delta} f\big]^q$  cannot be integrable in $\R^\N$. Therefore, we also need $\delta<0$.
\v
{\bf Case 2}~~Let $p>1$. Consider $f(x)=\chi_Q(x)|x|^{-\delta\big({p\over p-1}\big)}$. The norm inequality in (\ref{norm ineq L^1}) implies
\bel{characteristic again}
\sup_{Q\subset\R^\N}~|Q|^{{\alpha\over \N}-{1\over p}+{1\over q}}\left\{{1\over |Q|}\int_{Q} \left({1\over |x|}\right)^{\gamma q} dx\right\}^{1\over q}\left\{{1\over |Q|}\int_{Q} \left({1\over |x|}\right)^{\delta \big({p\over p-1}\big)} dx\right\}^{p-1\over p}~<~\infty.
\eeq
We essentially need $\gamma<\N/q$ and $\delta<\N\left({p-1\over p}\right)$ for which $|x|^{-\gamma q}$ and $|x|^{-\delta\big({p\over p-1}\big)}$ are locally integrable. By changing dilations inside (\ref{characteristic again}), we find ${\alpha\over\N}={1\over p}-{1\over q}+{\gamma+\delta\over\N}$ is an necessity.
Moreover, we claim ${\alpha\over \N}\ge {1\over p}-{1\over q}$. Together with ${\alpha\over\N}={1\over p}-{1\over q}+{\gamma+\delta\over\N}$, we must have $\gamma+\delta\ge0$. 

Suppose ${\alpha\over \N}<{1\over p}-{1\over q}$. Let $Q$ centered on some $x_o\neq0$. By shrinking $Q$ to $x_o$ and applying Lebesgue's Differentiation Theorem, we find 
\[ 
\left\{{1\over |Q|}\int_{Q} \left({1\over |x|}\right)^{\gamma q} dx\right\}^{1\over q}\left\{{1\over |Q|}\int_{Q} \left({1\over |x|}\right)^{\delta\big({p\over p-1}\big)} dx\right\}^{p-1\over p}~=~|x_o|^{-(\gamma+\delta)}~>~0.
\]
On the other hand, 
$|Q|^{{\alpha\over \N}-{1\over p}+{1\over q}}\mt\infty$.
We reach a contradiction to (\ref{characteristic again}).

\subsection{Constraints in (\ref{constraints+formula}) imply the norm inequality in (\ref{norm ineq L^1})}
Consider
\bel{split}
\I_{\alpha\gamma\delta} f(x)~=~\U_1f(x)+\U_2f(x)+\U_3f(x),\qquad f\ge0
\eeq 
where
\bel{U_1}
\begin{array}{cc}\ds
\U_1 f(x)~=~\int_{\R^\N} f(y)\Omega_1(x,y)dy,
\\\\ \ds
\Omega_1(x,y)~=~\left\{\begin{array}{lr}\ds\left({1\over|x|}\right)^\gamma\left({1\over|x-y|}\right)^{\N-\alpha}\left({1\over|y|}\right)^\delta,\qquad\hbox{$0\leq{|y|\over |x|}\leq{1\over 2}$},
\\ \ds ~~~~~~~~~~~~~~~~~~0,~~~~~~~~~~~~~~~~~~~~~~~~~~~~~~~~~\hbox{${|y|\over |x|}>{1\over 2}$};\end{array}\right.
\end{array}
\eeq

\bel{U_2}
\begin{array}{cc}\ds
\U_2 f(x)~=~\int_{\R^\N} f(y)\Omega_2(x,y)dy,
\\\\ \ds
\Omega_2(x,y)~=~\left\{\begin{array}{lr}\ds\left({1\over|x|}\right)^\gamma\left({1\over|x-y|}\right)^{\N-\alpha}\left({1\over|y|}\right)^\delta,\qquad~~\hbox{${|y|\over |x|}\ge2$},
\\ \ds ~~~~~~~~~~~~~~~~~~0,~~~~~~~~~~~~~~~~~~~~~~~~~~~~~\hbox{$0\leq{|y|\over |x|}<2$}\end{array}\right.
\end{array}
\eeq
and
\bel{U_3}
\begin{array}{cc}\ds
\U_3 f(x)~=~\int_{\R^\N} f(y)\Omega_3(x,y)dy,
\\\\ \ds
\Omega_3(x,y)~=~\left\{\begin{array}{lr}\ds\left({1\over|x|}\right)^\gamma\left({1\over|x-y|}\right)^{\N-\alpha}\left({1\over|y|}\right)^\delta,\qquad~~~~~~~~~~\hbox{${1\over 2}<{|y|\over |x|}<2$},
\\ \ds ~~~~~~~~~~~~~~~~~~0,~~~~~~~~~~~~~~~~~~~~~~~~~~~~~~~\hbox{$0\leq{|y|\over |x|}\leq{1\over 2}$ ~or~ ${|y|\over |x|}\ge{1\over 2}$}.
\end{array}\right.
\end{array}
\eeq
We aim to show that each $\U_i f, i=1,2,3$ satisfies the $\L^p\mt\L^q$-norm inequality in (\ref{norm ineq L^1}). The first two can be obtained by refining the proof of Stein and Weiss \cite{Stein-Weiss} and using the two lemmas from the previous section. The crucial part of our estimates comes when dealing with $\U_3 f$ for which we go through an interpolation argument of changing measures.
\v

{\bf Case One}~~Let $1\leq p=q<\infty$. The homogeneity condition ${\alpha\over\N}={1\over p}-{1\over q}+{\gamma+\delta\over\N}$ implies $\alpha=\gamma+\delta$.

 Note that $|x-y|\geq{1\over 2}|x|$ if ${|y|\over|x|}\leq{1\over 2}$. From (\ref{U_1}), we find 
\bel{K1 estimate}
\Omega_1(x,y)
~\leq~2^{\N-\alpha}\left\{\begin{array}{lr}\ds \left({1\over|x|}\right)^{\N-\alpha+\gamma} \left({1\over|y|}\right)^\delta,\qquad~~~~~\hbox{\small{$0\leq{|y|\over |x|}\leq{1\over 2}$}},
\\ \ds ~~~~~~~~~~~~~0,~~~~~~~~~~~~~~~~~~~~~~~~~~~~\hbox{${|y|\over |x|}>{1\over 2}$}.
\end{array}\right.
\eeq
Because $\delta<\N\left({p-1\over p}\right)$, we have
\bel{J1}
\begin{array}{lr}\ds
\A_1~\doteq~\int_0^{\infty}\Omega_1(1,t)t^{\N\left({p-1\over p}\right)-1}dt
~\leq~2^{\N-\alpha}\int_0^{1\over2}t^{\N\left({p-1\over p}\right)-\delta-1}dt~<~\infty.
\end{array}
\eeq
By applying {\bf Lemma One}, we obtain
\bel{U1 a Ineq}
\left\|\U_1f\right\|_{\L^p(\R^\N)}~\leq~\C_{\alpha~\delta~p}~\left\|f\right\|_{\L^p(\R^\N)},\qquad 1\leq p<\infty.
\eeq
\vsk
On the other hand,  $|x-y|\geq{1\over 2}|y|$ if ${|y|\over |x|}\geq2$. From (\ref{U_2}), we find 
\bel{K2 estimate}
\Omega_2(x,y)
~\leq~2^{\N-\alpha}\left\{\begin{array}{lr}\ds
\left({1\over|x|}\right)^{\gamma} \left({1\over|y|}\right)^{\N-\alpha+\delta},\qquad~~\hbox{${|y|\over |x|}\ge2$},
\\ \ds ~~~~~~~~~~0,~~~~~~~~~~~~~~~~~~~~~~~~\hbox{$0\leq{|y|\over |x|}<2$}
\end{array}\right.
\eeq
Because $\gamma<\N/q=\N/p$ and $\alpha=\gamma+\delta$, we have 
\bel{J2}
\begin{array}{lr}\ds
\A_2~\doteq~\int_{0}^\infty \Omega_2(1,t)t^{\N\left({p-1\over p}\right)-1}dt
\\\\ \ds~~~~
~\leq~2^{\N-\alpha}\int_2^\infty t^{\N\left({p-1\over p}\right)-\N+\alpha-\delta-1}dt~=~2^{\N-\alpha}\int_2^\infty t^{-{\N\over p}+\gamma-1}dt
~<~\infty.
\end{array}
\eeq
By applying {\bf Lemma One}, we obtain
\bel{U2 a Ineq}
\left\|\U_2f\right\|_{\L^p(\R^\N)}~\leq~\C_{\alpha~\gamma~p}~\left\|f\right\|_{\L^p(\R^\N)},\qquad 1\leq p<\infty.
\eeq

Let $\N\ge2$. Write $x=R\xi$ and $y=r\eta$ for $\xi,\eta\in\mathds{S}^{\N-1}$. Recall $\Omega_3$ defined in (\ref{U_3}). We have
\bel{U3 estimate}
\begin{array}{lr}\ds
\U_3f(x)~=~\int_{\R^\N} \Omega_3(x,y)f(y)dy
\\\\ \ds
~~~~~~~~~~~~=~\int_{\mathds{S}^{\N-1}}\int_0^\infty \Omega_3(R\xi, r\eta)f(r\eta)r^{\N-1}drd\sigma(\eta)
\\\\ \ds
~~~~~~~~~~~~=~\int_{\mathds{S}^{\N-1}}\int_0^\infty \Omega_3(\xi, t\eta)f(tR\eta)t^{\N-1}dtd\sigma(\eta)
\\ \ds~~~~~~~~~~~~~~~~~~~~~~~~~~~~~
\hbox{\small{$r=tR$, $\Omega_3$ is homogeneous of degree $-\N$}}
\\\\ \ds
~~~~~~~~~~~~=~\int_{\mathds{S}^{\N-1}}\int_{1\over 2}^2 {1\over |\xi-t\eta|^{\N-\alpha}}f(tR\eta)t^{\N-1-\delta}dtd\sigma(\eta)\qquad\hbox{\small{by (\ref{U_3})}}
\\\\ \ds
~~~~~~~~~~~~=~\int_{\mathds{S}^{\N-1}}\int_{1\over 2}^2 {1\over \big[(\xi-t\eta)\cdot (\xi-t\eta)\big]^{\N-\alpha\over2}}f(tR\eta)t^{\N-1-\delta}dtd\sigma(\eta)
\\\\ \ds
~~~~~~~~~~~~=~\int_{\mathds{S}^{\N-1}}\int_{1\over 2}^2 {1\over |1-2t\xi\cdot\eta+t^2|^{\N-\alpha\over2}}f(tR\eta)t^{\N-1-\delta}dtd\sigma(\eta)
\\\\ \ds
~~~~~~~~~~~~\lesssim~\int_{\mathds{S}^{\N-1}}\int_{1\over 2}^2{1\over\Delta^{\N-\alpha}(t,\xi,\eta)}f(tR\eta) dt d\sigma(\eta).
\end{array}
\eeq
For $\N=1$, take $d\sigma$ to be point measure on $1$ and $-1$ inside (\ref{U3 estimate}). We find
\bel{U3 estimate N=1}
\begin{array}{lr}\ds
\U_3f(x)~\lesssim~\int_{1\over 2}^2 |1-t|^{\alpha-1} \Big[ f(tR)+f(-tR)\Big]dt.
\end{array}
\eeq

From (\ref{U3 estimate}), we have
\bel{U3 a Ineq}
\begin{array}{lr}\ds
\left\|\U_3f\right\|_{\L^p(\R^\N)}~\lesssim~\left\{\int_{\mathds{S}^{\N-1}}\int_0^\infty\left\{\int_{\mathds{S}^{\N-1}}\int_{1\over 2}^2{1\over\Delta^{\N-\alpha}(t,\xi,\eta)}f(tR\eta)dt d\sigma(\eta)\right\}^pR^{\N-1}dRd\sigma(\xi)\right\}^{1\over p}
\\\\ \ds
\lesssim\int_{1\over 2}^2\left\{\int_0^\infty\int_{\mathds{S}^{\N-1}}\left\{\int_{\mathds{S}^{\N-1}}{1\over\Delta^{\N-\alpha}(t,\xi,\eta)}f(tR\eta) d\sigma(\eta)\right\}^pd\sigma(\xi)  R^{\N-1}dR\right\}^{1\over p} dt
\\ \ds
~~~~~~~~~~~~~~~~~~~~~~~~~~~~~~~~~~~~~~~~~~~~~~~~~~~~~~~~~~~~~~~~~~\hbox{\small{by Minkowski integral inequality}}
\\\\ \ds
\lesssim\int_{1\over 2}^2\left\{\int_0^\infty\int_{\mathds{S}^{\N-1}}\left\{\int_{\mathds{S}^{\N-1}}{\big[f(tR\eta)\big]^p\over\Delta^{\N-\alpha}(t,\xi,\eta)} d\sigma(\eta)\right\}\left\{ \int_{\mathds{S}^{\N-1}} {1\over\Delta^{\N-\alpha}(t,\xi,\eta)} d\sigma(\eta)\right\}^{p-1}d\sigma(\xi)  R^{\N-1}dR\right\}^{1\over p} dt
\\ \ds~~~~~~~~~~~~~~~~~~~~~~~~~~~~~~~~~~~~~~~~~~~~~~~~~~~~~~~~~~~~~~~~~~~~~~~~~~~~~~~~~~~~~~~~~~~~~~~~~~~~~~~~~~~~~~~~~\hbox{\small{by H\"{o}lder inequality}}
\\\\ \ds
\lesssim\int_{1\over 2}^2\left\{\int_0^\infty\int_{\mathds{S}^{\N-1}}\left\{\int_{\mathds{S}^{\N-1}}{\big[f(tR\eta)\big]^p\over\Delta^{\N-\alpha}(t,\xi,\eta)} d\sigma(\eta)\right\}|1-t|^{\left[{\alpha-\N\over\N}\right](p-1)}d\sigma(\xi)  R^{\N-1}dR\right\}^{1\over p} dt
\\ \ds~~~~~~~~~~~~~~~~~~~~~~~~~~~~~~~~~~~~~~~~~~~~~~~~~~~~~~~~~~~~~~~~~~~~~~~~~~~~~~~~~~~~~~~~~~~~~~~~~~~~~~~~~\hbox{\small{by {\bf Lemma Two}}}
\\\\ \ds
=\int_{1\over 2}^2\left\{\int_0^\infty\int_{\mathds{S}^{\N-1}}\left\{|1-t|^{\left[{\alpha-\N\over\N}\right](p-1)}\int_{\mathds{S}^{\N-1}}{1\over\Delta^{\N-\alpha}(t,\xi,\eta)} d\sigma(\xi) \right\}  \big[f(tR\eta)\big]^p  R^{\N-1}dRd\sigma(\eta)\right\}^{1\over p} dt
\\\\ \ds
\lesssim\int_{1\over 2}^2\left\{|1-t|^{\left[{\alpha-\N\over\N}\right]p}\int_0^\infty \int_{\mathds{S}^{\N-1}} \big[f(tR\eta)\big]^p  R^{\N-1}dRd\sigma(\eta)\right\}^{1\over p} dt\qquad\hbox{\small{by {\bf Lemma Two}}}
\\\\ \ds
=\left\{\int_0^\infty \int_{\mathds{S}^{\N-1}} \big[f(tR\eta)\big]^p d\sigma(\eta) R^{\N-1}dR\right\}^{1\over p} \int_{1\over 2}^2 |1-t|^{\alpha-\N\over\N}dt
\\\\ \ds
= \left\| f\right\|_{\L^p(\R^\N)} \int_{1\over 2}^2 |1-t|^{\alpha-\N\over\N}dt
~\leq~\C_\alpha \left\| f\right\|_{\L^p(\R^\N)},\qquad 1\leq p<\infty.
\end{array}
\eeq
Moreover, by using (\ref{U3 estimate N=1}), we find
\bel{U3 a Ineq N=1}
\begin{array}{lr}\ds
\left\|\U_3f\right\|_{\L^p(\R)}~\lesssim~\left\{\int_0^\infty\left\{ \int_{1\over 2}^2 |1-t|^{\alpha-1} \Big[ f(tR)+f(-tR)\Big]dt\right\}^p dR\right\}^{1\over p}
\\\\ \ds~~~~~~~~~~~~~~~~~~
~\leq~\int_{1\over 2}^2\left\{\int_0^\infty \Big[ f(tR)+f(-tR)\Big]^p dR\right\}^p  |1-t|^{\alpha-1}dt
\\ \ds~~~~~~~~~~~~~~~~~~~~~~~~~~~~~~~~~~~~~~~~~~~\hbox{\small{by Minkowski integral inequality}}
\\\\ \ds~~~~~~~~~~~~~~~~~~
~=~\left\| f\right\|_{\L^p(\R)} \int_{1\over 2}^2 |1-t|^{\alpha-1}dt
\\\\ \ds~~~~~~~~~~~~~~~~~~
~\leq~\C_\alpha \left\| f\right\|_{\L^p(\R)},\qquad 1\leq p<\infty.
\end{array}
\eeq

{\bf Case Two}~~ Consider $1\leq p<q<\infty$. Assert
\bel{V delta}
\V_\delta f(x)~=~|x|^{-\N+\delta}\int_{|y|<|x|}|y|^{-\delta}f(y)dy,\qquad\hbox{\small{$\delta<\N\left({p-1\over p}\right)$}}.
\eeq
We claim
\bel{a}
\left\|\V_{\delta}f\right\|_{\L^p(\R^\N)}~\leq~\C_{\delta~p}~\|f\|_{\L^p(\R^\N)},\qquad  1\leq p<\infty.
\eeq
Write
\bel{V Omega}
\V_\delta f(x)~=~\int_{\R^n}\Omega(|x|,|y|)f(y)dy
,\qquad
\Omega(u,v)~=~\left\{
\begin{array}{lr}\ds
u^{-\N+\delta}v^{-\delta}\qquad \hbox{if}~~v<u
\\ \ds
~~~~~~0\qquad~~~~~~~~\hbox{otherwise}.
\end{array}
\right.
\eeq
Observe that $\Omega$ in (\ref{V Omega}) is homogeneous of degree $-\N$. Moreover, 
\bel{V Omega est}
\int_0^\infty \Omega(1,t) t^{\N\big({p-1\over p}\big)-1}dt~=~\int_0^1 t^{\N\big({p-1\over p}\big)-\delta-1} dt~<~\infty 
\eeq
provided by $\delta<\N\left({p-1\over p}\right)$. {\bf Lemma One}  implies (\ref{a}).

On the other hand, $\V_\delta f$ defined in (\ref{V delta}) satisfies 
\bel{V est norm}
\begin{array}{lr}\ds
\V_\delta f(x)~\leq~|x|^{-\N+\delta} \left\{\int_{|y|<|x|} |y|^{-\delta \big({p\over p-1}\big)} dy\right\}^{p-1\over p} \left\| f\right\|_{\L^p(\R^\N)}\qquad \hbox{\small{by H\"{o}lder inequality}}
\\\\ \ds~~~~~~~~~~~
~\leq~\C_{\delta~p}~|x|^{-{\N\over p}}  \left\| f\right\|_{\L^p(\R^\N)},\qquad 1\leq p<\infty.
\end{array}
\eeq
Recall $\U_1 f$ defined in (\ref{U_1}). Note that $0\leq\frac{|y|}{|x|}\leq\frac{1}{2}$ implies $\frac{1}{2}|x|\leq|x|-|y|\leq|x-y|$.  Let $f,g\ge0$ and $f\in\L^p(\R^\N)$, $g\in\L^{q\over q-1}(\R^\N)$. We have
\bel{U_1 bilinear}
\begin{array}{lr}\ds
\int_{\R^\N} \U_1 f(x) g(x)dx~=~
\int_{\R^\N}\left\{\int_{|y|\leq\frac{1}{2}|x|}{f(y)g(x)\over~|x|^{\gamma}~|x-y|^{\N-\alpha}~|y|^{\delta}~}dy\right\} dx
\\\\ \ds~~~~~~~~~~~~~~~~~~~~~~~~~~~~~
~\lesssim~\int_{\R^\N}\left\{\int_{|y|<|x|}{f(y)g(x)\over~|x|^{\gamma+\N-\alpha}~|y|^{\delta}~}dy\right\}dx
\\\\ \ds~~~~~~~~~~~~~~~~~~~~~~~~~~~~~
~=~\int_{\R^\N}|x|^{\alpha-(\gamma+\delta)}g(x)\left\{|x|^{-\N+\delta}\int_{|y|<|x|}f(y)|y|^{-\delta}dy\right\}dx
\\\\ \ds~~~~~~~~~~~~~~~~~~~~~~~~~~~~~
~=~\int_{\R^\N}|x|^{\alpha-(\gamma+\delta)}g(x)\V_\delta f(x)dx
\\\\ \ds~~~~~~~~~~~~~~~~~~~~~~~~~~~~~
~\leq~\left\{\int_{\R^\N}~|x|^{\big[\alpha-(\gamma+\delta)\big]q}\Big(\V_{\delta}f\Big)^q(x)dx\right\}^{1\over q}\left\|g\right\|_{\L^{q\over q-1}(\R^\N)}~~\hbox{\small{by H\"{o}lder inequality}}.
\end{array}
\eeq

Let ${\alpha\over \N}={1\over p}-{1\over q}+{\gamma+\delta\over\N}, 1\leq p<q<\infty$.
We find
\bel{U_1 est q>1}
\begin{array}{lr}\ds
\left\{\int_{\R^\N}~|x|^{\big[\alpha-(\gamma+\delta)\big]q}\Big(\V_{\delta}f\Big)^q(x)dx\right\}^{1\over q}
\\\\ \ds
\leq~\left\{\int_{\R^\N}~|x|^{\big[\alpha-(\gamma+\delta)\big]q}|x|^{-\N\big[{q\over p}-1\big]}\left\| f\right\|_{\L^p(\R^\N)}^{q-p}\Big(\V_{\delta}f\Big)^p(x)dx\right\}^{1\over q}\qquad\hbox{\small{by (\ref{V est norm})}}
\\\\ \ds
=~\left\| f\right\|_{\L^p(\R^\N)}^{1-{p\over q}} \left\{\int_{\R^\N} \Big(\V_\delta f\Big)^p(x)dx\right\}^{1\over q}
~\leq~\C_{\delta~p}~\left\| f\right\|_{\L^p(\R^\N)}\qquad\hbox{\small{by (\ref{a})}}.
\end{array}
\eeq
From (\ref{U_1 bilinear})-(\ref{U_1 est q>1}), we conclude
\bel{U_1 Result}
\left\|\U_1 f\right\|_{\L^q(\R^\N)}~\leq~\left\| f\right\|_{\L^p(\R^\N)},\qquad 1\leq p<q<\infty.
\eeq
Consider
\bel{V gamma}
\V_\gamma g(x)~=~|x|^{-\N+\gamma}\int_{|y|<|x|}|y|^{-\gamma}g(y)dy,\qquad\hbox{\small{$\gamma<{\N\over q}$}}.
\eeq
We claim
\bel{b}
\left\|\V_\gamma g\right\|_{\L^{q\over q-1}(\R^\N)}~\leq~\C_{\gamma~q}~\|g\|_{\L^{q\over q-1}(\R^\N)},\qquad  1< q<\infty.
\eeq
Write
\bel{V Omega again}
\V_\gamma g(x)~=~\int_{\R^n}\Omega(|x|,|y|)g(y)dy
,\qquad
\Omega(u,v)~=~\left\{
\begin{array}{lr}\ds
u^{-\N+\gamma}v^{-\gamma}\qquad \hbox{if}~~v<u
\\ \ds
~~~~~~0\qquad~~~~~~~~\hbox{otherwise}.
\end{array}
\right.
\eeq
Observe that $\Omega$ in (\ref{V Omega again}) is homogeneous of degree $-\N$. Moreover, 
\bel{V Omega est again}
\int_0^\infty \Omega(1,t) t^{{\N\over q}-1}dt~=~\int_0^1 t^{{\N\over q}-\gamma-1} dt~<~\infty 
\eeq
provided by $\gamma<{\N\over q}$. {\bf Lemma One}  implies (\ref{b}).

On the other hand, $\V_\gamma g$ defined in (\ref{V gamma}) satisfies
\bel{V gamma est}
\begin{array}{lr}\ds
\V_\gamma g(y)~=~|y|^{-\N+\gamma}\int_{|x|<|y|}|x|^{-\gamma}g(x)dx
\\\\ \ds~~~~~~~~~~~
~\leq~|y|^{-\N+\gamma} \left\{\int_{|x|<|y|} |x|^{-\gamma q}dx\right\}^{1\over q}\left\| g\right\|_{\L^{q\over q-1}(\R^\N)}\qquad\hbox{\small{by H\"{o}lder inequality}}
\\\\ \ds~~~~~~~~~~~
~\leq~\C_{\gamma~q}~|y|^{-\N\big({q-1\over q}\big)}\left\| g\right\|_{\L^{q\over q-1}(\R^\N)}. 
\end{array}
\eeq

Recall $\U_2 f$ defined in (\ref{U_2}). Note that $\frac{|y|}{|x|}\geq2$ implies $\frac{1}{2}|y|\leq|y|-|x|\leq|x-y|$. We have
\bel{U_2 bilinear}
\begin{array}{lr}\ds
\int_{\R^\N} \U_2 f(x) g(x)dx~=~\int_{\R^\N}\left\{ \int_{|y|\ge2|x|}{f(y)\over~|x|^{\gamma}~|x-y|^{\N-\alpha}|y|^\delta}dy\right\}g(x) dx
\\\\ \ds~~~~~~~~~~~~~~~~~~~~~~~~~~~~~~
~\lesssim~\int_{\R^\N}\left\{\int_{|y|>|x|}{f(y)g(x)\over|x|^{\gamma}~|y|^{\N-\alpha+\delta}}dy\right\}dx.
\end{array}
\eeq
By using (\ref{U_2 bilinear}) and Tonelli's theorem, we have
\bel{U_2 estimate}
\begin{array}{lr}\ds
\int_{\R^\N} \U_2 f(x) g(x)dx
~\lesssim~\int_{\R^\N}|y|^{\alpha-(\gamma+\delta)}f(y)\left\{|y|^{-\N+\gamma}\int_{|x|<|y|}g(x)|x|^{-\gamma}dx\right\}dy
\\\\ \ds~~~~~~~~~~~~~~~~~~~~~~~~~~~~~~
~=~\int_{\R^\N}|y|^{\alpha-(\gamma+\delta)}f(y)\V_\gamma g(y)dy
\\\\ \ds~~~~~~~~~~~~~~~~~~~~~~~~~~~~~~
~\leq~\left\| f\right\|_{\L^p(\R^\N)} \left\{\int_{\R^\N} |y|^{\big[\alpha-(\gamma+\delta)\big]{p\over p-1}} \Big(\V_\gamma g\Big)^{p\over p-1}(y)dy\right\}^{p-1\over p}\qquad\hbox{\small{by H\"{o}lder inequality}}.
\end{array}
\eeq
Let ${\alpha\over \N}={1\over p}-{1\over q}+{\gamma+\delta\over\N}, 1\leq p<q<\infty$.
For $p=1$, we find
\bel{U_2 p=1 est}
\begin{array}{lr}\ds
\left\| |y|^{\alpha-(\gamma+\delta)} \V_\gamma g(y)\right\|_{\L^\infty(\R^\N)}~\leq~\C_{\gamma~q} ~|y|^{\alpha-(\gamma+\delta)} |y|^{-\N\big({q-1\over q}\big)}\left\| g\right\|_{\L^{q\over q-1}(\R^\N)}\qquad\hbox{\small{by (\ref{V gamma est})}}
\\\\ \ds~~~~~~~~~~~~~~~~~~~~~~~~~~~~~~~~~~~~~~~~~
~=~\C_{\gamma~q}~\left\| g\right\|_{\L^{q\over q-1}(\R^\N)}.
\end{array}
\eeq
For $p>1$, write
\[\begin{array}{lr}\ds
\left\{\int_{\R^\N} |y|^{\big[\alpha-(\gamma+\delta)\big]{p\over p-1}} \Big(\V_\gamma g\Big)^{p\over p-1}(y)dy\right\}^{p-1\over p}
\\\\ \ds
=~\left\{\int_{\R^\N} |y|^{\big[\alpha-(\gamma+\delta)\big]{p\over p-1}} \Big(\V_\gamma g\Big)^{\big[{p\over p-1}-{q\over q-1}\big]}(y)\Big(\V_\gamma g\Big)^{q\over q-1}(y)dy\right\}^{p-1\over p}.
\end{array}\]
By using (\ref{V gamma est}) again, we have
\bel{U_2 est}
\begin{array}{lr}\ds
|y|^{\big[\alpha-(\gamma+\delta)\big]{p\over p-1}} \Big(\V_\gamma g\Big)^{\big[{p\over p-1}-{q\over q-1}\big]}(y)
~\leq~\C_{\gamma~q}~|y|^{\big[\alpha-(\gamma+\delta)\big]{p\over p-1}}|y|^{-\N\big({q-1\over q}\big)\big[{p\over p-1}-{q\over q-1}\big]}\left\| g\right\|_{\L^{q\over q-1}(\R^\N)}^{{p\over p-1}-{q\over q-1}}
\\\\ \ds~~~~~~~~~~~~~~~~~~~~~~~~~~~~~~~~~~~~~~~~~~~~~~~~~~~
~=~\C_{\gamma~q}~|y|^{\N\big[{1\over p}-{1\over q}\big]{p\over p-1}-\N\big({q-1\over q}\big)\big[{p\over p-1}-{q\over q-1}\big]} \left\| g\right\|_{\L^{q\over q-1}(\R^\N)}^{{p\over p-1}-{q\over q-1}}
\\\\ \ds~~~~~~~~~~~~~~~~~~~~~~~~~~~~~~~~~~~~~~~~~~~~~~~~~~~
~=~\C_{\gamma~q} \left\| g\right\|_{\L^{q\over q-1}(\R^\N)}^{{p\over p-1}-{q\over q-1}}.
\end{array}
\eeq
From (\ref{U_2 estimate})-(\ref{U_2 est}), together with the $\L^{q\over q-1}$-estimate in (\ref{b}), we conclude
\bel{U_2 Result}
\left\|\U_2 f\right\|_{\L^q(\R^\N)}~\leq~\left\| f\right\|_{\L^p(\R^\N)},\qquad 1\leq p<q<\infty.
\eeq
Recall $\U_3 f$ defined in (\ref{U_3}). We have
\bel{U_3 dominate}
\begin{array}{lr}\ds
\U_3 f(x)~=~\int_{{1\over 2}|x|<|y|<2|x|} f(y) \left({1\over |x|}\right)^\gamma\left({1\over |x-y|}\right)^{\N-\alpha} \left({1\over |y|}\right)^\delta dy
\\\\ \ds~~~~~~~~~~~
~\lesssim~\int_{\R^\N} f(y) \left({1\over |x-y|}\right)^{\N-(\alpha-\gamma-\delta)} dy.
\end{array}
\eeq
Note that ${\alpha-\gamma-\delta\over\N}={1\over p}-{1\over q}, 1\leq p<q<\infty$. Define
\bel{I_alpha-gamma-delta}
\I_{\alpha-\gamma-\delta} f(x)~=~\int_{\R^\N} f(y) \left({1\over |x-y|}\right)^{\N-(\alpha-\gamma-\delta)} dy.
\eeq
For $p>1$, we have
\bel{I_abc Est}
\left\| \I_{\alpha-\gamma-\delta} f\right\|_{\L^q(\R^\N)} ~\leq~\C_{p~q}~\left\| f\right\|_{\L^p(\R^\N)}.
\eeq
This is a classical result due to Hardy, Littlewood-Sobolev \cite{Hardy-Littlewood} and Sobolev \cite{Sobolev}.
\begin{remark} For $p=1$, we have $\I_{\alpha-\gamma-\delta}\colon \L^1(\R^\N)\mt\L^{q,\infty}(\R^\N)$. See chapter V of Stein \cite{Stein}.
\end{remark}
Given $E\subset\R^\N$, denote $\vol\{ E\}=\int_E dx$.
From (\ref{U_3 dominate})-(\ref{I_alpha-gamma-delta}) and {\bf Remark 3.1}, we have
\bel{weak type est}
\begin{array}{lr}\ds
\lambda^q \vol\left\{ x\in\R^\N\colon \U_3 f(x)>\lambda\right\}~\leq~\lambda^q \vol\left\{ x\in\R^\N\colon \I_{\alpha-\gamma-\delta} f(x)>\lambda\right\}
\\\\ \ds~~~~~~~~~~~~~~~~~~~~~~~~~~~~~~~~~~~~~~~~~~~~~~~
~\lesssim~\left\| f\right\|_{\L^1(\R^\N)}^q,\qquad \lambda>0.
\end{array}
\eeq
By replacing $f(x)$ with $f(x)|x|^\delta$ inside (\ref{weak type est}), we obtain
\bel{weak type Est}
\lambda ~\vol\left\{ x\in\R^\N\colon \int_{{1\over 2}|x|<|y|<2|x|} f(y) \left({1\over |x|}\right)^\gamma\left({1\over |x-y|}\right)^{\N-\alpha}  dy>\lambda\right\}^{1\over q}~\lesssim~\int_{\R^\N} f(x)|x|^\delta dx
\eeq
for every $\lambda>0$.

Recall $\delta<0$. Let $\delta_1<\delta<\delta_2<0$ of which $\delta_i, i=1,2$ are close to $\delta$. We find
\bel{delta constraints i=1,2}
{\alpha\over \N}~=~1-{1\over q_i}+{\gamma+\delta_i\over\N},\qquad i~=~1,2
\eeq
for some $q_1>q>q_2>1$. 

By carrying out the same argument as (\ref{U_3 dominate})-(\ref{weak type Est}), we simultaneously have
\bel{weak type Est i=1,2}
\begin{array}{lr}\ds
\lambda ~\vol\left\{ x\in\R^\N\colon \int_{{1\over 2}|x|<|y|<2|x|} f(y) \left({1\over |x|}\right)^\gamma\left({1\over |x-y|}\right)^{\N-\alpha}  dy>\lambda\right\}^{1\over q_i}~\lesssim~\int_{\R^\N} f(x)|x|^{\delta_i} dx,
\\ \ds~~~~~~~~~~~~~~~~~~~~~~~~~~~~~~~~~~~~~~~~~~~~~~~~~~~~~~~~~~~~~~~~~~~~~~~~~~~~~~~~~~~~~~~~~~~~~~~~~~~~~~~~~~~~~~~~~~~~~~~
 i~=~1,2
 \end{array}
\eeq
for every $\lambda>0$.

Next, we need to apply a Marcinkiewicz interpolation theorem of changing measures, due to Stein and Weiss \cite{S-W}. 

Let $\mu_i, i=1,2$ be two absolutely continuous measures satisfying
\bel{mu}
\mu_i(E)~=~\int_{E}~|x|^{\delta_i}dx,\qquad i~=~1,2.
\eeq
Define
\bel{mu_t}
\mu_t(E)=\int_{E}|x|^{\delta_1(1-t)}|x|^{\delta_2t}dx,\qquad \frac{1}{q_t}~=~\frac{1-t}{q_1}+\frac{t}{q_2},\qquad 0~\leq ~t~\leq~1.
\eeq
\v

{\bf Theorem B: ~~Stein and Weiss,~ 1958} \\
{\it 
Let $\T$ be a sub-linear operator, having the following properties:

 (1)~~The domain of $\T$ includes $\L^1(\R^\N, d\mu_1) \cap\L^1(\R^\N, d\mu_2)$.

 (2)~~If $f\in\L^1(\R^\N, d\mu_i)$, $i=1,2$, we have
\bel{weak type conditions}
\lambda~\vol\left\{x\in\R^\N:\left|\T f(x)\right|>\lambda\right\}^{1\over q_i}~\leq~\C\int_{\R^\N} |f(x)| d\mu_i(x),\qquad i=1,2.
\eeq
Then, 
\bel{T result}
\left\|\T f\right\|_{\L^{q_t}(\R^\N)}~\leq~\C\int_{\R^\N} |f(x)| d\mu_t(x),\qquad 0<t<1.
\eeq}
\begin{remark} {\bf Theorem B} in the more general setting of measurable spaces can be found in the paper by Stein and Weiss \cite{S-W}.
\end{remark}
Recall ${\alpha\over \N}=1-{1\over q}+{\gamma+\delta\over\N}$ and (\ref{delta constraints i=1,2}).
There is a $0<t<1$ such that
$\delta=(1-t)\delta_1+t\delta_2$ and ${1\over q}={1-t\over q_1}+{t\over q_2}$.
By using (\ref{weak type Est i=1,2}) and applying {\bf Theorem B}, we obtain
\bel{U_3 Result}
\left\|\U_3 f\right\|_{\L^q(\R^\N)}~\leq~\left\| f\right\|_{\L^1(\R^\N)}.
\eeq

\section{Proof of Theorem Two}
\setcounter{equation}{0}
\subsection{Norm inequality in (\ref{Norm Ineq L^1}) implies constraints in (\ref{Formula})-(\ref{Constraints})}
Let $\II_{\alpha\gamma\delta}$ defined in (\ref{I_alpha f}) for $0<\alpha_i<\N_i, i=1,2,\ldots,n$ and $\gamma,\delta\in\R$.  First, it is clear that $\gamma<{\N\over q}, \delta<0, \gamma+\delta\ge0, {\alpha\over\N}=1-{1\over q}+{\gamma+\delta\over\N}$  are still required because $\left({1\over |x-y|}\right)^{\N-\alpha}\leq \prod_{i=1}^n\left({1\over |x_i-y_i|}\right)^{\N_i-\alpha_i}$. Next, we show $\alpha_i-\N_i<\delta, i=1,2,\ldots,n$ as an extra necessary condition. 

Let $\Q\doteq \Q_1\times\Q_2\times\cdots\times\Q_n$ where $\Q_i\subset\R^{\N_i}, i=1,2,\ldots,n$ is a cube. We write $\Q=\Q_i\times\Q_i'$ where $\Q_i'=\bigotimes_{j\neq i} \Q_j$. Denote $x=(x_i,x'_i)\in\R^{\N_i}\times\R^{\N-\N_i}$. 

Let $f$  be an indicator function supported in $\Q$. 
The $\L^1\mt\L^q$-norm inequality in (\ref{Norm Ineq L^1}) implies
\bel{Characteristic}
\begin{array}{lr}\ds
\A_q^{\alpha~\gamma~\delta}(\Q)~=~\prod_{i=1}^n|\Q_i|^{{\alpha_i\over \N_i}-1+{1\over q}}\left\{{1\over |\Q|}\int_{\Q} \left({1\over |x|}\right)^{\gamma q} dx\right\}^{1\over q}\left\{{1\over |\Q|}\int_{\Q} \left({1\over |x|}\right)^\delta dx\right\}
\\\\ \ds~~~~~~~~~~~~~~~~
~\doteq~\prod_{i=1}^n|\Q_i|^{{\alpha_i\over \N_i}-1+{1\over q}}~\B_q^{\alpha~\gamma~\delta}(\Q)~\leq~\C_{\alpha~\gamma~\delta~q},\qquad \Q\subset\R^\N.
\end{array}
\eeq
We claim
${\alpha_i\over \N_i}-1+{1\over q}\ge0, i=1,2,\ldots,n$.
Suppose ${\alpha_i\over \N_i}-1+{1\over q}<0$ for some $i\in\{1,2,\ldots,n\}$. Consider $\Q$ centered on the origin of $\R^\N$. Let $|\Q_i|^{1\over \N_i}=\lambda$, $0<\lambda<1$ and $|\Q_j|^{1\over\N_j}=1$ for $j\neq i$. 
By shrinking $\Q_i$ to $0\in\R^{\N_i}$ and applying Lebesgue's Differentiation Theorem, we find
\bel{int ter prod}
\begin{array}{lr}\ds
\lim_{\lambda\mt0}\B_q^{\alpha~\gamma~\delta}(\Q)~=~
\left\{{1\over |\Q'_i|}\int_{\Q'_i} \left({1\over |x'_i|}\right)^{\gamma q} dx'_i\right\}^{1\over q}\left\{{1\over |\Q'_i|}\int_{\Q'_i} \left({1\over |x'_i|}\right)^\delta dx'_i\right\}~>~0.
\end{array}
\eeq
Consequently, $\A_q^{\alpha~\gamma~\delta}(\Q)$ in (\ref{Characteristic}) diverges to infinity as $\lambda\mt0$.

Let
\bel{Q^k}
\Q_i^k~=~\Q_i\cap\Big\{ 2^{-k-1}\leq|x_i|<2^{-k}\Big\},\qquad k\ge0,\qquad i=1,2,\ldots,n.
\eeq 
Given $i\in\{1,2,\ldots,n\}$, we assert $|\Q_i|^{1\over \N_i}=1$ and $|\Q_j|^{1\over \N_j}=\lambda$ for  $j\neq i$. Write
\bel{Chara rewrite}
\begin{array}{lr}\ds
\prod_{i=1}^n|\Q_i|^{q\big[{\alpha_i\over \N_i}-1+{1\over q}\big]}\left\{{1\over |\Q|}\int_{\Q} \left({1\over |x|}\right)^{\gamma q} dx\right\}\left\{{1\over |\Q|}\int_{\Q} \left({1\over |x|}\right)^\delta dx\right\}^q
~=~
\\\\ \ds
\sum_{k\ge0}\prod_{j\neq i} |\Q_j|^{q\big[{\alpha_j\over \N_j}-1+{1\over q}\big]} \left\{\prod_{j\neq i} {1\over |\Q_j|}\iint_{\Q_i^k\times\bigotimes_{j\neq i}\Q_j} {\sqrt{|x_i|^2+\sum_{j\neq i}|x_j|^2}}^{-\gamma q} dx_i \prod_{j\neq i} dx_j\right\}
\left\{{1\over |\Q|}\int_{\Q} \left({1\over |x|}\right)^\delta dx\right\}^q
\\\\ \ds
~\doteq~\sum_{k\ge0} \A_k(\lambda).
\end{array}
\eeq
By applying Lebesgue's Differentiation Theorem, we have
\bel{Chara EST}
\begin{array}{lr}\ds
\lim_{\lambda\mt0}~ \prod_{j\neq i} {1\over |\Q_j|}\iint_{\Q_i^k\times\bigotimes_{j\neq i}\Q_j} {\sqrt{|x_i|^2+\sum_{j\neq i}|x_j|^2}}^{-\gamma q} dx_i \prod_{j\neq i} dx_j
~=~ \int_{\Q_i^k} \left({1\over |x_i|}\right)^{\gamma q} dx_i.
\end{array}
\eeq
Suppose that there is a $ j\neq i$ such that ${\alpha_j\over \N_j}-1+{1\over q}>0$. 
We find $\A_k(0)=0$ for every $k\ge0$. Moreover, each $\A_k(\lambda)$ is H\"{o}lder continuous for $\lambda\ge0$ whose exponent is strict positive depending on ${\alpha_j\over \N_j}-1+{1\over q}$. Recall the inequality in (\ref{Characteristic}).  
 For every $\lambda>0$, $\sum_{k\ge0}\A_k(\lambda)\leq\C_{\alpha~\gamma~\delta~q}$. Consequently, $\sum_{k\ge0}\A_k(\lambda)$ is continuous at $\lambda=0$. We have 
\bel{Chara Lim}
\lim_{\lambda\mt0} ~\sum_{k\ge0}\A_k(\lambda)~=~0.
\eeq
A direct computation shows
\bel{Characteristic Est Case1}
\begin{array}{lr}\ds
\prod_{i=1}^n |\Q_i|^{q\big[{\alpha_i\over \N_i}-1+{1\over q}\big]}\left\{{1\over|\Q|}\int_{\Q} \left({1\over|x|}\right)^{\gamma q} dx\right\}\left\{{1\over |\Q|}\int_{\Q} \left({1\over |x|}\right)^\delta dx\right\}^q
\\\\ \ds
~\ge~\C_{\delta~q}\prod_{j\neq i}\lambda^{q\big[\alpha_j-\N_j+{\N_j\over q}\big]}\int_{\Q_i}{\sqrt{|x_i|^2+\sum_{j\neq i}\lambda^2}}^{-\gamma q}dx_i
\\\\ \ds
~\ge~\C_{\delta~q}\prod_{j\neq i}\lambda^{q\big[\alpha_j-\N_j+{\N_j\over q}\big]}\int_{0<|x_i|\leq\lambda}\left({1\over \lambda}\right)^{\gamma q} dx_i
~=~\C_{\gamma~\delta~q}~\lambda^{\N_i-\gamma q+\sum_{j\neq i} q\big[\alpha_j-\N_j+{\N_j\over q}\big]}.
\end{array}
\eeq
From (\ref{Chara Lim})-(\ref{Characteristic Est Case1}), by using ${\alpha\over\N}=1-{1\over q}+{\gamma+\delta\over\N}$, we find
\bel{Constraint est.1}
\begin{array}{lr}\ds
{\N_i\over q}-\gamma+\sum_{j\neq i} \alpha_j-\N_j+{\N_j\over q}~>~0\qquad\Longrightarrow\qquad \alpha_i~<~{\N\over q}-\gamma+\alpha-\N+\N_i~=~\N_i+\delta.
\end{array}
\eeq
On the other hand, suppose ${\alpha_j\over \N_j}-1+{1\over q}=0$ for every $j\neq i$. Similar to (\ref{Characteristic Est Case1}), we find
\bel{Characteristic Est Case2}
\begin{array}{lr}\ds
\prod_{i=1}^n |\Q_i|^{q\big[{\alpha_i\over \N_i}-1+{1\over q}\big]}\left\{{1\over|\Q|}\int_{\Q} \left({1\over|x|}\right)^{\gamma q} dx\right\}\left\{{1\over |\Q|}\int_{\Q} \left({1\over |x|}\right)^\delta dx\right\}^q
\\\\ \ds
~\ge~\C_{\delta~q}\prod_{j\neq i}\lambda^{q\big[\alpha_j-\N_j+{\N_j\over q}\big]}\int_{\Q_i}{\sqrt{|x_i|^2+\sum_{j\neq i}\lambda^2}}^{-\gamma q}dx_i
\\\\ \ds
~\ge~\C_{\delta~q}\int_{\lambda<|x_i|\leq1}\left({1\over |x_i|}\right)^{\gamma q} dx_i.
\end{array}
\eeq
The last integral in  (\ref{Characteristic Est Case2}) converges as $\lambda\mt0$. We must have $\gamma q<\N_i$. By using ${\alpha\over\N}=1-{1\over q}+{\gamma+\delta\over\N}$ and taking into account  ${\alpha_j\over \N_j}-1+{1\over q}=0$ for every $j\neq i$, we find
\bel{Constraint est.2}
\begin{array}{lr}\ds
\alpha_i~=~\N_i-{\N_i\over q}+\gamma+\delta\qquad\Longrightarrow\qquad \alpha_i~<~\N_i+\delta.
\end{array}
\eeq
\subsection{Constraints in (\ref{Formula})-(\ref{Constraints}) imply the norm inequality in (\ref{Norm Ineq L^1})}
Now, we prove the $\L^1\mt\L^q$-norm inequality in (\ref{Norm Ineq L^1}) for $\alpha,\gamma,\delta$ satisfying (\ref{Formula})-(\ref{Constraints}).
By symmetry, we assume $|x_1-y_1|\ge|x_i-y_i|, i=2,\ldots,n$.
Denote $\t$ to be an ($n-1$)-tuple: 
$\left(2^{-t_2},\dots,2^{-t_n}\right)$ where $t_i\ge0, ~i=2,\ldots,n$. Let $f\ge0$. We define
\bel{Delta II}
\Delta_\t\II_{\alpha\gamma\delta} f(x)~=~\int_{\Gamma_\t(x)}f(y)\left({1\over |x|}\right)^\gamma\prod_{i=1}^n\left(\frac{1}{|x_i-y_i|}\right)^{\N_i-\alpha_i}\left({1\over |y|}\right)^\delta dy,\qquad x\neq0
\eeq
where
\bel{cone}
\Gamma_\t(x)~=~\left\{y\in\R^\N:~2^{-t_i-1}<\frac{|x_i-y_i|}{|x_1-y_1|}\leq2^{-t_i},~i=2,\ldots,n\right\}.
\eeq
\begin{figure}[h]
\centering
\includegraphics[scale=0.25]{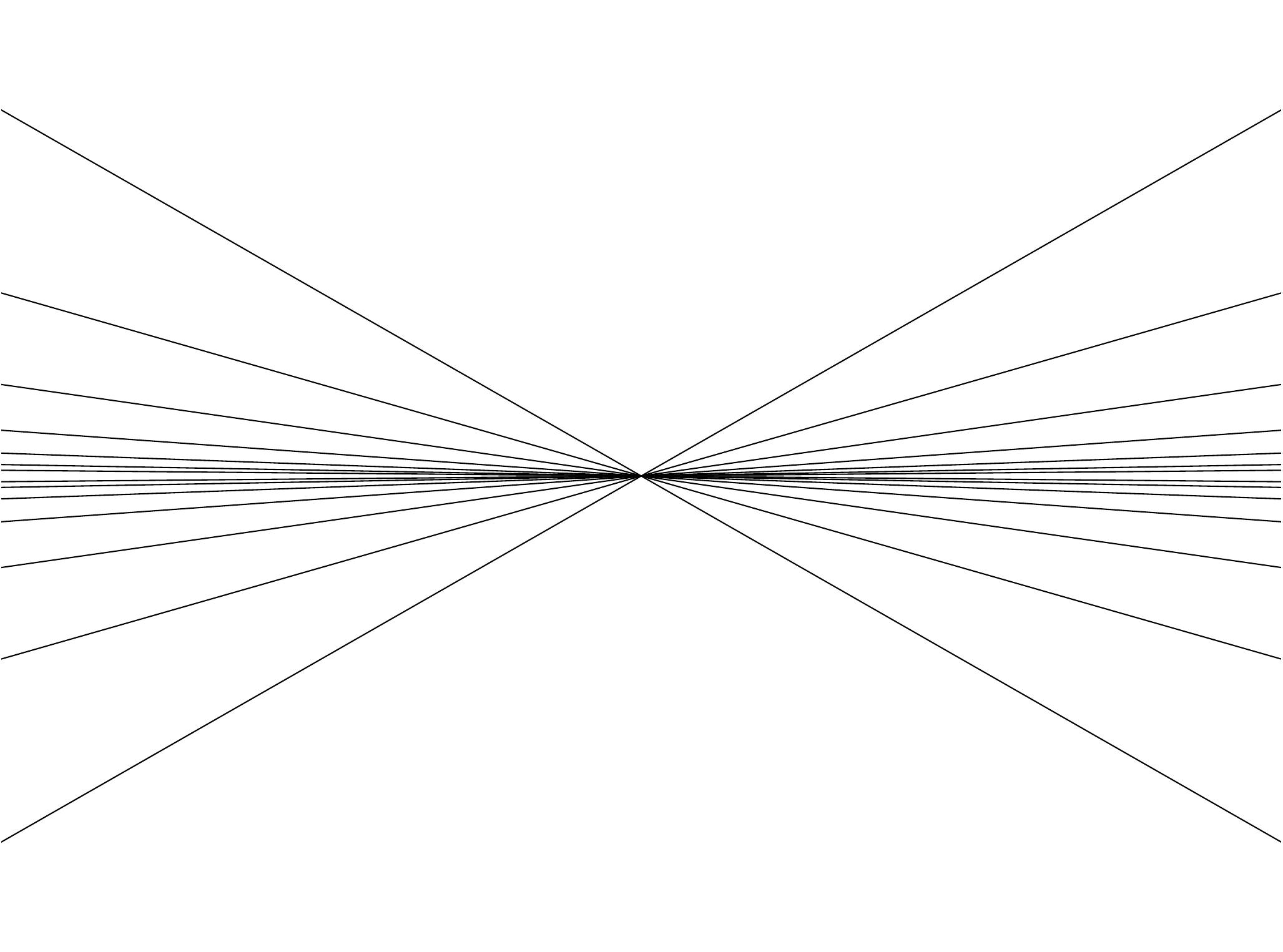}
\end{figure}

In particular, we write $\Gamma_o(x)=\Gamma_\t(x)$ for $t_2=\cdots=t_n=0$.

Let
\bel{n-dila}
\t x~=~\left(x_1,2^{-t_2}x_2,\dots, 2^{-t_n}x_n\right),\qquad \t^{-1} x~=~\left(x_1,2^{t_2}x_2,\dots, 2^{t_n}x_n\right).
\eeq
We have
\bel{EST.1 cone}
\begin{array}{lr}\ds
\Delta_\t\II_{\alpha\gamma\delta} f(x)~=~\int_{\Gamma_\t(x)}f(y)\left(\frac{1}{|x|}\right)^{\gamma}\prod_{i=1}^n\left(\frac{1}{|x_i-y_i|}\right)^{\N_i-\alpha_i}\left(\frac{1}{|y|}\right)^\delta dy
\\\\ \ds
~=~\int_{\Gamma_o(x)}f(\t^{-1} y)\left(\frac{1}{|\t^{-1} x|}\right)^{\gamma}\prod_{i=1}^n\left(\frac{1}{|2^{t_i}x_i-2^{t_i}y_i|}\right)^{\N_i-\alpha_i}\left(\frac{1}{|\t^{-1} y|}\right)^\delta2^{\N_2t_2+\cdots+\N_nt_n}dy
\\ \ds~~~~~~~~~~~~~~~~~~~~~~~~~~~~~~~~~~~~~~~~~~~~~~~~~~~~~~~~~~~~~~~~~~~~~~~~~~~~~~~~~~~~~~~~~~~x\rightarrow\t^{-1} x, y\rightarrow\t^{-1} y
\\\\ \ds
~=~2^{\alpha_2t_2+\cdots+\alpha_nt_n}\int_{\Gamma_o(x)}f(\t^{-1} y)\left(\frac{1}{|\t^{-1} x|}\right)^{\gamma}\prod_{i=1}^n\left(\frac{1}{|x_i-y_i|}\right)^{\N_i-\alpha_i}\left(\frac{1}{|\t^{-1} y|}\right)^\delta dy.
\end{array}
\eeq
Denote
\bel{t_mu}
t_\mu~\doteq~\max\left\{t_i: i=2,\dots,n\right\}.
\eeq
Recall $\gamma+\delta\ge0$ and $\delta<0$. Because $t_i\geq0$, $i=2,\dots,n$, we find
\bel{weight compara}
\left|\t^{-1} x\right|^\gamma~\geq~\left|x\right|^\gamma,\qquad \left|\t^{-1} y\right|^\delta~\geq~2^{t_\mu\delta}\left|y\right|^\delta.
\eeq
From (\ref{EST.1 cone}) to (\ref{weight compara}), we further have
\bel{EST.2 cone}
\begin{array}{lr}\ds
\Delta_\t\II_{\alpha\gamma\delta} f(x)~\leq~2^{\alpha_1t_1+\alpha_2t_2+\cdots+\alpha_nt_n}2^{-t_\mu\delta}\int_{\Gamma_o(x)}f(\t^{-1} y)\prod_{i=1}^n\left(\frac{1}{|x|}\right)^{\gamma}\left(\frac{1}{|x_i-y_i|}\right)^{\N_i-\alpha_i}\left(\frac{1}{|y|}\right)^\delta dy
\\\\ \ds~~~~~~~~~~~~~~~~~~
~\lesssim~2^{\alpha_1t_1+\alpha_2t_2+\cdots+\alpha_nt_n}2^{-t_\mu\delta}\int_{\R^\N}f(\t^{-1} y)\left(\frac{1}{|x|}\right)^{\gamma}\left(\frac{1}{|x-y|}\right)^{\N-\alpha}\left(\frac{1}{|y|}\right)^\delta dy.
\end{array}
\eeq 
By using (\ref{EST.2 cone}) and applying {\bf Theorem One}, we obtain
\bel{Delta Norm Ineq}
\begin{array}{lr}\ds
\left\|\Delta_\t\II_{\alpha\gamma\delta} f\right\|_{\L^q(\R^\N)}~\leq~\C_{\alpha~\gamma~\delta~q}~2^{\alpha_2t_2+\cdots+\alpha_nt_n}2^{-t_\mu\delta}\int_{\R^\N}  f(\t^{-1} y) dy
\\\\ \ds~~~~~~~~~~~~~~~~~~~~~~~~~~~~
~=~\C_{\alpha~\gamma~\delta~q}~2^{(\alpha_2-\N_2)t_2+\cdots+(\alpha_n-\N_n)t_n}2^{-t_\mu\delta}
\left\| f\right\|_{\L^1(\R^\N)}.
\end{array}
\eeq
Recall $\alpha_i-\N_i<\delta, i=1,2,\ldots,n$. Let
$\ve=\min\left\{\N_\mu-\alpha_\mu+\delta, ~\N_i-\alpha_i,~i\neq\mu\right\}>0$.
From (\ref{Delta Norm Ineq}), we find
\bel{Delta norm decays}
\begin{array}{lr}\ds
\left\|\Delta_\t\II_{\alpha\gamma\delta} f\right\|_{\L^q(\R^\N)}~\leq~\C_{\alpha~\gamma~\delta~q}~2^{-\ve\sum_{i=2}^n t_i}~\left\|f\right\|_{\L^1(\R^\N)}.
\end{array}
\eeq
Lastly,  by applying Mikowski inequality, we finish the proof of {\bf Theorem Two}.

{\small School of Mathematical Sciences, Zhejiang University}\\
{\small email: 12063001@zju.edu.cn}

{\small Department of Mathematics, Westlake University}\\
{\small email: wangzipeng@westlake.edu.cn}

\end{document}